\documentclass[12pt]{article} 
\title{Some differential equations in SDG}
\author{Anders Kock and Gonzalo E. Reyes}
\newtheorem{theorem}{Theorem}
\newtheorem{proposition}[theorem]{Proposition}
\newtheorem{lemma}[theorem]{Lemma}

\newtheorem{corollary}[theorem]{Corollary}

\begin{document}
\maketitle
\section*{Introduction} We intend to comment on some of those aspects
of the theory of differential equations which we think are clarified
(for us, at least) by means of the synthetic method. By this, we
understand that the objects under consideration are seen as objects
in one sufficiently rich category (model for SDG), allowing us, for
instance, to work with nilpotent numbers, say $d\in R$ with $d^2 =0$;
but the setting should also permit the formation of function spaces,
so that some of the methods of functional analysis, become available,
in particular, the theory of distributions.

The specific topics we treat are generalities on vector fields and
the solutions of corresponding first- and second-order ordinary
differential equations; and also some partial differential equations,
which can be seen in this light, the wave- and heat-equation on some
simple spaces, like the line $R$. For these equations, distribution
theory is not just a tool, but is rather the essence of the matter,
since what develops through time, is a distribution (of heat, say),
which, as stressed by Lawvere,  is an {\em extensive} quantity, and
as such behaves covariantly, unlike density {\em functions} (which 
behave contravariantly); and the distributions may have no density 
function, in particular in the setting of model for SDG where all 
functions are {\em smooth}.

When we consider these partial differential equations, we shall
follow an
old practice and sometimes denote derivative $d/dt$ with respect
to ``time'' by a dot, $\dot{f}$, whereas differential operators with 
respect to space variables are denoted $\partial f/ \partial x$, $f'$, 
$\Delta (f)$, etc.

\medskip

We want to thank Henrik Stetk\ae r for useful conversations on the
topic of distributions.

\section{Generalities on actions}

Recall that an {\em action} of a set (object) $D$ on a set (object)
$M$ is a map $X :D\times M \to M$, and a {\em homomorphism} of actions 
$(M, X ) \to (N, Y)$ is a map $f: M \to N$ with $f(X (d , m))= Y(d , 
f(m))$ for all $m \in M$ and $d\in D$.

\medskip
The category of actions by a set $D$ form a topos; we shall in
particular be interested in the exponent formation in this topos, when
the action in the exponent is {\em invertible}. An action $X :D\times
M \to M$ is called {\em invertible}, if for each $d \in D$, $X (d,
-): M\to M$ is invertible. In this case, the exponent $(N,Y)^{(M, X
)}$ may be described as $N^M$ equipped with the following action by
$D$: an element $d\in D$ acts on $\beta :M \to N$ by ``conjugation":
$$\beta \mapsto Y_d \circ \beta \circ (X _d)^{-1},$$
where $Y_d$ denotes $Y(d, -):N\to N$, and similarly for $X _d$.

In the applications below, $D$ is the usual set of square zero
elements in $R$.
It is a {\em pointed} object, pointed by $0\in D$, and the actions
$X :D\times M \to M$  we consider, are {\em pointed} actions in the
sense that $X (0,m)=m$ for all $m\in M$, or equivalently,
$X _0 :M \to M$ is the identity map on $M$.
A pointed action, in this situation, is the same thing as a vector
field on $M$, cf.\ \cite{catdyn}.

\medskip
In the above situation, if $X $ and $Y$ are pointed actions, then so
is the exponent described. The  pointed actions likewise form a
topos, and the exponent described is then also the exponent in the
category of pointed actions; cf.\  \cite{kr2}.

For the case of vector fields seen as actions by $D$, we want to
describe the ``streamlines" generated by a vector field in abstract
action-theoretic terms; this is going to involve a certain
``universal" action $(\tilde{R}, \Delta )$: $\tilde{R}$ is an
``infinitesimally open subset" of $R,$ i.e., whenever $x\in \tilde{R}$
then $x+d\in \tilde{R}$ for every $d\in D.$ The main examples of such
subsets are $R$ itself, the non-negative numbers $R_{\geq 0},$ open
intervals, and the set $D_{\infty }$ of all nilpotent elements of the
number line. The action $\Delta$ is the vector field $\partial /
\partial x$, meaning the map $D\times \tilde{R} \to \tilde{R}$ given
by $(d,t)\mapsto d+t$.  (So it is not to be confused with the Laplace 
operatot $\Delta$, to be considered later.) The main property to be assumed is 
that the individual $\Delta _d$'s are homomorphisms of $D$-actions 
(which is a commutativity requirement); the structure of $\tilde{R}$ 
could probably be derived from this, but we shall be content with {\em 
assuming} that $\tilde{R}$ is an additively written monoid, and that 
$D\subseteq \tilde{R}$ (with the $0$ of $D$ also being the zero of the 
monoid).

\medskip
First, if $(M, X) $ is a set with an action, a homomorphism $f:
(\tilde{R}, \Delta ) \to (M, X )$ is to be thought of as a particular solution
of the differential equation given by $X$, with initial value $f(0)$,
or as a ``streamline" for the vector field $X$, starting in $f(0)$.
One wants, however, also to include dependence on initial value into
the notion of solution, and so one is led to consider maps
$$F: \tilde{R}\times M \to M,$$
satisfying at least $F(d,m)=X (d,m)$ for all $d\in D$ and $m\in M$;
we shall consider and compare the following further conditions
(universally quantified over all $d\in D$, $t, s\in \tilde{R}$, $m\in
M$):

\begin{equation}F(\Delta (d,t),m))= X (d, F(t,m));
\label{one} \end{equation}
this is the main one, the two following conditions are included for
systematic reasons only:

\begin{equation}F(\Delta (d,t),m))=F(t, X (d,m)),
\label{two} \end{equation}

\begin{equation}F(t, X (d,m))=X (d, F(t,m))
\label{three} \end{equation}
Finally, one may consider the following equation
\begin{equation}F(t+s,m)= F(t, F(s,m)).
\label{four} \end{equation}
Writing $X_d$ for the map $X(d, -):M\to M$, and similarly for $F$,
condition (\ref{one}) may be rewritten as

$$F_{\Delta (d, t)}=X_d \circ F_t$$
The others may be rewritten in a similar way. For instance
(\ref{four}) may be rewritten as
$$F_{t+s}= F_t \circ F_s $$
Equation (\ref{one}) expresses that, for each fixed $m\in M$, the map
$F(-,m):\tilde{R} \to M$ is a homomorphism (and thus, by virtue of
$F(0,m)=m$, a ``solution with initial value $m$"). Writing the action
of $D$ in terms of the symbol $\cdot$, we may write it
$F(d\cdot t , m) = d\cdot F(t,m).$ Equation (\ref{two}) expresses a
certain bi-homogeneity condition of $F$, $F(d\cdot t, m )= F(t,
d\cdot m)$ ; (\ref{three}) says that for fixed $t\in R$, $F(t,-):M
\to M$ is an endomorphism of $D$-actions, $F(t,d\cdot m) = d\cdot
F(t,m)$. Finally (\ref{four}) is the usual condition for action af a
monoid on a set $M$. Clearly, it implies all the others.

\medskip
Let $X$ be a vector field on $M$, thought of as a first-order
differential equation.  We say that the map $F:\tilde{R}\times M\to M$ is a 
{\it complete solution} or simply a {\it solution} if $F_d=X_d$ and 
$F$ satisfies (\ref{one}).  A solution in this sense does not satisfy 
the other conditions (\ref{two})-(\ref{four}), but it does, provided 
that $M$ satisfies a certain axiom (reflecting, synthetically, 
validity of the uniqueness assertion for solutions of differential 
equations on $M$).  --- The axiom in question is the following

\medskip
\noindent {\bf Uniqueness property for $M$}:

\noindent {\it If $X$ is a $D$-action on $M$, and $f, g: \tilde{R}\to M$
are homomorphisms of actions, with $f(0) = g(0)$, then $f=g$.}

\medskip

Note that the validity of the axiom, for a given $M$, depends on the
choice of $\tilde{R}, \Delta$. For instance, we shall prove below
that it holds for any microlinear $M$ if $\tilde{R}$ is taken to be
$D_{\infty }$ (and $\Delta =\partial /\partial x $).

\begin{proposition} Let $X$ be a vector field on $M$ and assume that
$M$ satisfies the uniqueness axiom. Then any solution
$F:\tilde{R}\times M \to M$ of the differential equation $X$
satisfies properties (\ref{two}) and (\ref{three}). Furthermore, if
$\tilde{R}$ is a monoid (under +) then $F$ also satisfies
(\ref{four}).
\end{proposition}

{\bf Proof.} Since the proofs are quite similar, we shall do only
(\ref{four}).
Fix $m\in M$ and $s\in \tilde{R}$ and define the couple of functions
$f,\;g:\tilde{R}\to M$ by the formulas

$$\left\{ \begin{array}{ll}
f(t)=F_{t+s}(m)\\
g(t)=F_t\circ F_s(m)
\end{array}
\right.$$
We have to check that $f$ and $g$ are homomorphisms of $D$-actions,
i.e., they satisfy (\ref{one}). Let us do this for the first

$$\begin{array}{llll}
f(t+d) & = & F_{(t+d)+s}(m) \\
       & = & F_{d+(t+s)}(m) \\
       & = & F_d\circ F_{t+s}(m) \\
       & = & X_d\circ f(t).
\end{array}
$$
The proof that $g$ is a homomorphism is similar. Thus, the equality of the two 
expressions follows from the uniqueness property assumed for $M$.

\medskip

Recall that a vector field $X$ on $M$ is called {\em integrable} if
there exists a solution $F: \tilde{R}\times M \to M$. If we assume
the uniqueness property, the equation (\ref{four}) holds; if further
the commutative monoid structure $+$ on $\tilde{R}$ actually is a
group structure, then (\ref{four}) implies that the action is
invertible, with $X_{-d}$ as $X_d ^{-1}$ (in fact $F_{-d} = F_d
^{-1}$). Of course, both the uniqueness property and the question
whether or not the vector field $X$ is integrable, depends on which
$\tilde{R}$ is considered. In particular, we shall say that $X$ is
{\it formally integrable} or {\it has a formal solution} if $X$ is
integrable for $\tilde{R}=D_{\infty}$ (which is a group under
addition). For the case of $M=R^n$, this amounts to integration by
formal power series, whence the terminology.

\begin{theorem}\label{formal} 
The uniqueness property holds for any microlinear object, (for  
$\tilde{R}=D_{\infty}$).  Furthermore, every vector field on a 
microlinear object is formally integrable.  Thus, every vector field 
on a microlinear object has a unique formal solution.
\end{theorem}
{\bf Proof.} We need to recall some infinitesimal objects from the 
literature on SDG, cf.\ e.g.\ \cite{lav}.  Besides $D\subseteq R$, 
consisting of $d\in R$ with $d^2 =0$, we have $D^n \subseteq R^n$, the 
$n$-fold product of $D$ with itself.  It has the subobject 
$D(n)\subseteq D^n$ consisting of those $n$-tuples $(d_1 , \ldots ,d_n 
)$ where $d_i\cdot d_j =0$ for all $i, j$.  There is also the object 
$D_n \subseteq R$ consisting of $\delta \in R$ with $\delta ^{n+1} 
=0$; $D_{\infty }$ is the union of all the $D_n$'s.  If $(d_1 , \ldots 
,d_n )\in D^n$, then $d_1 + \ldots + d_n \in D_n$.

--- Now, let $M$ be a microlinear object, and $X$ a vector field on it.  We 
first recall that if $d_1 ,d_2 \in D$ have the property that $d_1 + 
d_2 \in D$, then $X_{d_1}\circ X_{d_2} = X_{d_1 + d_2 }$.  (For 
microlinear objects perceive $D(2)$ to be a pushout over $\{0\}$ of 
the two inclusions $D\to D(2)$, and clearly both expressions given 
agree if either $d_1 =0$ or $d_2 =0$.) In particular, $X_{d_1}$ and 
$X_{d_2}$ commute.  But more generally, \begin{lemma} If $X$ is a 
vector field on a microlinear object and $d_1 , d_2 \in D$, the maps 
$X_{d_1}$ and $X_{d_2}$ commute.
\label{lemm3}\end{lemma}
{\bf Proof.} This is a consequence of the theory of Lie brackets, cf.\
e.g.\  \cite{lav} 3.2.2, namely $[X,X] = 0$.

\medskip

Likewise
\begin{lemma} If $X$ is a vector field on a microlinear object and
$d_1,\dots, d_n \in D$ are such that $d_1 +\dots+ d_n =0$, then

$$X_{d_1}\circ\dots \circ X_{d_n} = 1_M$$
(= the identity map on $M$). In particular,
$(X_d)^{-1} = X_{-d}.$
\label{lemm4}
\end{lemma}
{\bf Proof.}
We first prove that $R$, and hence any microlinear
object, perceives $D_n$ to be the orbit space of $D^n $ under the
action of the symmetric group ${\bf S}_n$ in $n$ letters: Assume that $p:
D^n \to R$ coequalizes the action, i.e. is symmetric in the $n$
arguments. By the basic axiom of SDG, $p$ may be written in the form
$$p(d_1 ,\dots, d_n ) = \sum_{Q\subseteq \{1,\dots, n\}} a_Q d^Q$$
for unique $a_Q$'s in $R$ (where $d^Q$ denotes $\prod _{i\in Q} d_i$).
We claim that $a_Q = a_{\pi (Q)}$ for every $\pi \in {\bf S}_n$.
Indeed,

$$\sum _Q a_Q d^Q = \pi (\sum _Q a_Q d^Q )$$
since $p$ is symmetric. But
$$ \pi (\sum _Q a_Q d^Q ) = \sum a_Q d^{\pi (Q)}= \sum _Q a_{\pi ^{-1}
(Q)} d^{Q}.$$
By comparing coefficients and using uniqueness of coefficients, we
conclude $a_Q = a_{\pi (Q)}$, and this shows that $p$ is (the
restriction to $D^n$ of) a symmetric polynomial $R^n \to R$. By
Newton's theorem (which holds internally), $p$ is a polynomial in the
elementary symmetric polynomials $\sigma _i$. Recall that $\sigma _1
(d_1,\dots,d_n )= d_1 +\dots +d_n$: and each $\sigma _i$, when
restricted to $D^n$, is a function of $\sigma_1$, since $d_1 ^2 = 0$;
e.g.
$$\sigma _2 (d_1 , d_2 ) = \sum d_i d_j = \frac{1}{2} (d_1 +\dots +
d_n )^2 =
\frac{1}{2} (\sigma _1 (d_1,\dots,d_n ))^2.$$

Now consider, for fixed $m\in M$, the map $p: D^n \to M$ given by
\newline $(d_1 ,\dots, d_n ) \mapsto X_{d_1} \circ\dots  \circ
X_{d_n} (m)$. By Lemma \ref{lemm3}, this map is invariant under the
symmetric group ${\bf S}_n$ (recall that this group is generated by
transpositions), so there is a unique $\phi :D_n \to M$ such that
$$\phi (d_1 +\dots +d_n ) = X_{d_1}\circ\dots \circ X_{d_n} (m).$$
So if $d_1 +\dots +d_n  = 0$, $X_{d_1}\circ\dots \circ X_{d_n} (m) =
\phi (0) = \phi (0 +\dots + 0)= X_0 \circ\dots  \circ X_0 (m) = m$.
This proves the Lemma.

\medskip

We can now prove the Theorem. We need to define $F_t : M\to M$ when
$t\in D_{\infty}$. Assume for instance that $t\in D_n $. By
microlinearity of $M$, $M$ perceives $D_n$ to be the orbit space of
$D^n$ under the action of ${\bf S}_n$ (see the proof of Lemma
\ref{lemm4}), via the map $(d_1 ,\dots ,d_n )\mapsto d_1 +\dots
+d_n$, so we are forced to define $ F_t = X_{d_1} \circ\dots
X_{d_n}$ if $F$ is to extend $X$ and to satisfy (\ref{four}). The
fact that this is well defined independently of the choice of $n$ and
the choice of $d_1 ,\dots , d_n$ that add up to $t$ follows from
Lemma \ref{lemm4}.

\medskip
As a particular case of special importance, we consider a {\em
linear} vector field on a microlinear and Euclidean $R$-module $V$.  
To say that the vector field is linear is to say that its 
principal-part formation $V\to V$ is a linear map, $\Delta$, say.
We have then the following version of a classical result:

\begin{proposition} Let a linear vector field on a microlinear
Euclidean $R$-module $V$ be given by the linear map $\Delta : V\to V$. Then the
unique formal solution of the corresponding differential equation,
i.e., the equation $\dot{F}(t)=\Delta (F(t))$ with initial position $v$,
 is the map $D_{\infty} \times V \to V$ given by 
\begin{equation}(t,v)\mapsto e^{t\cdot \Delta }
(v),\label{exp1}\end{equation}
where the right hand side here means the sum of the following
``series" (which has only finitely many non-vanishing terms, since $t$
is assumed nilpotent):
$$v + t\Delta (v) +\frac{t^2}{2!}\Delta ^2 (v) + \frac{t^3}{3!}\Delta
^3 (v)+ \dots$$
\label{exp-sol}\end{proposition}
Here of course $\Delta ^2 (v)$ means $\Delta (\Delta (v))$, etc.

\medskip

{\bf Proof.} We have to prove that $\dot{F}(t)= \Delta (F(t))$. We
calculate the left hand side by differentiating  the
series term by term (there are only finitely many non-zero terms): 
$$\Delta (v) +\frac{2t}{2!}\cdot \Delta^2 (v) + \frac{3t^2}{3!} 
\Delta^3 (v) + ...  = \Delta (v+ t\cdot \Delta (v) +
\frac{t^2}{2!}\cdot \Delta ^2 (v) + ... )$$
using linearity of $\Delta$.  But this is just $\Delta$ applied to 
$F(t)$.

\medskip
There is an analogous result for second order differential  equations
of the form $\stackrel{\cdot \cdot}{F}(t)= \Delta (F(t))$ (with
$\Delta $ linear); the proof is similar and we omit it:

\begin{proposition} The formal solution of this second order
differential equation $\stackrel{\cdot \cdot}{F}= \Delta F$, with 
initial position $v$ and initial speed $w$, is given by $$F(t) = v + 
t\cdot w + \frac{t^2}{2!} \Delta (v) + \frac{t^3}{3!} \Delta (w) + 
\frac{t^4}{4!} \Delta ^2 (v) + \frac{t^5}{5!} \Delta ^2
(w)+ ... .$$
\label{exp-sol2}\end{proposition}

\section{Exponent vector fields}
In this section, we show that solutions of an exponent vector 
field may be obtained by conjugating solutions of the vector fields that 
make up the exponent. Furthermore, this method of conjugation is equivalent 
(under some conditions) to the method of change of variables, widely used 
to solve differential equations. 

\begin{theorem} Assume that $(M, X )$ and $(N,Y)$ are vector fields 
having solutions $F: \tilde{R} \times M \to M$ and $G: \tilde{R}
\times N \to N$, respectively, and assume that all $F_t$ are
invertible. Then a solution $H: \tilde{R}\times M\to M$ of the
exponent $(N,Y)^{(M, X )}$ is obtained as the map

$$H: \tilde{R} \times N^M \to N^M$$
given by conjugation: $H_t (\beta ) = G_t \circ \beta \circ F_t
^{-1}$.
\end{theorem}
{\bf Proof.} This is purely formal. For $\beta \in N^M$, we have
$$\begin{array}{llll}
(Y^{X})_d (H_t (\beta )) & = & Y_d \circ H_t (\beta ) \circ X _d
^{-1} \\
                         & = & Y_d \circ G_t \circ \beta \circ F_t
^{-1} \circ X _d ^{-1} \\

                         & = & G_{d+t} \circ \beta \circ F_{d+t}^{-1}
\\
                         & = & H_{d+t} (\beta ),
\end{array}$$ 
where in the third step we used the equation (\ref{one}) for $G$ and $F$, 
in the form $$G_{d+t} = Y_d \circ G_t,\; \mbox{ respectively } F_{d+t} 
= X _d \circ
F_t  ,$$
together with invertibility of $F_s$ for all $s$ and invertibility of
$X _d$.

A similar argument gives that if each of (\ref{two})-(\ref{four})
holds for both $F$ and $G$, then the corresponding property holds for
$H$.

\medskip
In most applications, the invertibility of the $F_t$ will be secured
by subtraction on $\tilde{R}$, with $F_t ^{-1} = F_{-t}$.

\medskip
Recall that an $R$-module $V$ is called {\em Euclidean} if the
canonical map $\alpha : V \times V \to V^D $ given by $\alpha
(u,v)(d)= u+d\cdot v$ is invertible; the composite of $\alpha ^{-1}$
with projection to the second factor,  $V^D \to V\times V \to V$ is
called {\em principal part} formation. If $X: V \to V^D$ is a vector
field on a Euclidean module $V$, we may compose it with principal
part formation to get a (not necessarily linear) map $\xi : V\to V$,
called the principal part of the vector field $X$; it is thus
characterized by the formula
$$X(v)(d)= v + d\cdot \xi (v).$$

Recall also that if $\beta : M\to V$ is any map into a Euclidean
$R$-module, and $X$ is a vector field on $M $, then the {\em
directional derivative} $D_X (\beta )$ of $\beta $ along $X$ is the
composite
$$M \stackrel{X}{\rightarrow} M^D \stackrel{\beta ^D}{\rightarrow}V^D
\rightarrow V,$$
where the last map is principal part formation. Using function
theoretic notation, $D_X (\beta )$ is characterized by validity of
the equation
$$\beta ( X(m,d )) = \beta (m) + d\cdot D_X (\beta ) (m),$$
for all $d\in D$, $m\in M$.

When $M$ itself is a Euclidean $R$ module, and $X$ has
principal part $\xi$, we usually write $D_{\xi} (\beta )$ instead of
$D_X (\beta )$.

\begin{proposition}Assume that $X_1$, $X_2$ are vector fields on
$M_1$, $M_2$, respectively, and that $H: M_1 \to M_2$ is a
homomorphism (i.e., it preserves the $D$-action). Let $V$ be a
Euclidean $R$-module. Then for any $u:M_2 \to V$,
$$D_{X_1} (u\circ H) = D_{X_2}(u)\circ H.$$
\label{lemma1}\end{proposition}
{\bf Proof.} This is a straightforward computation:
$$u(X_2 (H(m),d)) = u(H(m))+d\cdot D_{X_2}(H(m));$$
on the other hand
$$u(X_2 (H(m),d)) = u (H(X_1 (m,d))) = u(H(m))+ d\cdot D_{X_1}
(u\circ H)(m).$$
By comparing these two expressions we obtain the conclusion of the
Proposition.

\medskip
For any object $N$, let us consider its  ``zero vector field" $Z$ ,
i.e., $Z_d$ is the identity map on $N$, for all $d$. For a vector
field $X$ on an object $M$, we then also have the ``vertical" vector
field $Z\times X$ on $N\times M$.

If we have a complete solution $F: \tilde{R} \times M \to M$ of a
vector field $X$ on $M$, we may consider the map $\overline{F} :
\tilde{R} \times M  \to \tilde{R} \times M $ given by $(t,m) \mapsto
(t, F(t,m))$
\begin{proposition}
The map $\overline{F}$ thus described is an automorphism of the
vector field $Z\times X$ on $\tilde{R} \times M $.
\label{lemma2}\end{proposition}
{\bf Proof.} By a straightforward diagram chase, one sees that this
is a restatement of (\ref{three}).

\medskip
We now consider  solutions $F: \tilde{R}\times V  \to V$ for such
vector fields, so equation (\ref{one}) holds: $X_d \circ F_t =
F_{t+d}$. In terms of principal parts, this equation may be rewritten as

$$\dot{F}_t (v) = \xi (F_t (v)).$$
Similarly, equation (\ref{two}) may be written as
\begin{equation}\dot{F}_t = D_{\xi } (F_t
).\label{twoprime}\end{equation}

\medskip
Using directional derivatives, we can give a more familiar expression
to the vector field (1ODE) $Y^X$ considered above on the object
$N^M$, when the base $N$ is a microlinear Euclidean $R$-module $V$, 
and the exponent $M$ is mocrolinear. In fact, letting $\eta$ be the principal 
part of the vector field $Y$ on $N=V$, we have, for $u\in V^M$, $m\in 
M$, $d\in D$ (recall that $(X_d )^{-1}= X_{-d}$)

$$\begin{array}{llll}
(Y^X)_d (u)(m) & = & Y_d \circ u \circ X_{-d} (m) \\
               & = & u((X_{-d}(m)) + d\cdot \eta (u(X_{-d}(m))) \\
               & = & u(m) -d \cdot D_X (u)(m) + d \cdot \eta (u(m)) \\
               & = & u(m) + d\cdot [-D_X (u) (m) + \eta (u(m))] \\
\end{array}$$
(at the third equality sign, a cancellation of $d\cdot d$ took place in the 
last term)

\medskip
In other words, the principal part of $Y^X$ is $\theta : M\to V$
given by
$$\theta (m) = \eta (u(m)) - D_X (u) (m).$$
Recalling that the 1ODE corresponding to a vector field $X$
on a Euclidean $R$-module $V$ may be written as $\dot{x} = \xi (x)$ where 
$\xi$ is the principal part of $X$.  In these terms, the above equation may be 
rewritten (leaving out the $m$, and modulo some obvious abuse of 
notation) as
$$\dot{u} = \eta (u) - D_X (u),$$
or still, recalling that $\dot{(-)}$ is ``derivative with respect to
time",
$$\frac{\partial u}{\partial t} + D_X (u) = \eta (u).$$
This is a PDE of first order ``in time".

\medskip

The following may be seen as a generalization of (\ref{twoprime}),
and is a form of the chain rule. We consider a vector field $X$ on
$M$, with  solution $F: \tilde{R} \times M \to M$. Let $U: \tilde{R}
\times M  \to V$ be any function with values in a Euclidean
$R$-module.

\begin{proposition}
Under these circumstances, we have
 $$\frac{\partial}{\partial t} U(t, F_t (m))= \frac{\partial
U}{\partial t}(t, F_t (m)) + (D_{Z\times X} U )(t, F_t (m))$$
for all $t\in \tilde{R}$, $m\in M$.
\label{lemma3}\end{proposition}

{\bf Proof.} Since $F$ is a solution of $X$, $F_{t+d} = X_d \circ
F_t$, and so for any $t, t' \in \tilde{R}$ $(Z\times X)_d (t',F_t
(m)) = (t',F_{t+d}(m)$. Therefore, by definition of directional
derivative, $$U(t',F_{t+d}(m))= U(t',F_t (m))+ d\cdot (D_{Z\times
X}U)(t',F_t (m)).$$

 Putting $t'= t+d$, we thus have
$$U(t+d,F_{t+d}(m))= U(t+d,F_t (m))+ d\cdot (D_{Z\times X}U)(t+d,F_t
(m))$$
$$=U(t+d,F_t (m))+ d\cdot (D_{Z\times X}U)(t,F_t (m))$$
by a standard cancellation of two $d$'s, after Taylor expansion. Expanding the first term, we may continue:

$$=U(t, F_t (m)) + d\cdot \frac{\partial U}{\partial t}(t, F_t (m)) +
d\cdot (D_{Z\times X}U)(t,F_t (m)).$$
On the other hand,
$$U(t+d,F_{t+d}(m))= U(t, F_t (m)) + d\cdot \frac{\partial}{\partial
t} U(t, F_t (m));$$
comparing these two expressions gives the result.

\medskip
The method of change of variables has been used extensively to solve
differential equations. We shall prove that our method for solving
the exponential differential equation $Y^X$, where $X$ is an
integrable vector field on $M$, $Y$ an integrable vector field on a Euclidean
$R$-module, and where  $\tilde{R}$ is symmetric with respect to the origin (if
$t\in \tilde{R}$, then $-t\in \tilde{R}$), may be seen as an
application of the method of change of variables.  We let $\eta: V\to
V$ denote the principal part of $Y$, as before.  Let $F: \tilde{R}
\times M \to M$ be the assumed solution of $X$, and let $\overline{F}:
\tilde{R} \times M \to \tilde{R} \times M$ be the map
$$\overline{F}(t,m)= (t, F(-t,m))$$ 
Then $\overline{F}$ (which represents the change of variables $\tau =t$, $\mu = F(-t,m)$) is invertible.

\begin{theorem} (``Change of variables"). If $u: \tilde{R}\times M \to
V$ is a particular  solution of $Y^X$, or, equivalently, of
\begin{equation}\frac{\partial u}{\partial t} + D_X (u)= \eta(u),
\label{p.15}\end{equation}
then the unique map $U :\tilde{R}\times M \to V$ given as the
composite
$$\tilde{R}\times M
\stackrel{(\overline{F})^{-1}}{\to}\tilde{R}\times M
\stackrel{u}{\to} V$$
is a particular solution of $Y^Z$, or, equivalently, of
\begin{equation}\frac{\partial U}{\partial t} = \eta
(U),\label{p.15,2}\end{equation}
and vice versa.
\end{theorem}

{\bf Proof.} Since $u(t,m) = U(t, F_{-t}(m))$,
we have
$$\frac{\partial u}{\partial t} (t,m) = \frac{\partial}{\partial t}
U(t, F_{-t}(m)) = \frac{\partial U}{\partial t}(t, F_{-t}(m)) -
D_{Z\times X}U(t,F_{-t}(m)),$$
by the chain rule, Proposition \ref{lemma3}. On the other hand,
$\overline{F}$ is an automorphism of the vector field $Z\times X$, by
Proposition \ref{lemma2}, and so, by construction of $\overline{F}$
and Proposition \ref{lemma1},

$$D_{Z\times X} (u) = D_{Z\times X} (U\circ \overline{F}) =
(D_{Z\times X} U )\circ \overline{F}.$$
Therefore,
$$0= \frac{\partial u}{\partial t} + D_{Z\times X} (u) - g(u) $$
$$= \frac{\partial U}{\partial t}(t, \mu ) - D_{Z\times X} (U)(t, \mu
) + D_{Z\times X} (U)(t, \mu ) - \eta (U(t, \mu )),$$
where $\mu = F_{-t} (m)$, i.e., $U$ is solution of $$\frac{\partial
U}{\partial t} = \eta (U),$$
proving the theorem (the vice versa part follows because
$\overline{F}$ is invertible).

\medskip

{\bf Example.} Let $D$ be the set of elements of square zero in $R$,
as usual. It carries a vector field, namely  the map
$e: D\times D \to D$ given by $(d, \delta )\mapsto (1+d)\cdot
\delta$. It is easy to see that this vector field is integrable, with
complete solution
$E:R\times D \to D$ given by $(t,\delta )\mapsto e^t \cdot \delta $.
Now consider the tangent vector bundle $M^D$ on $M$. The zero vector
field $Z$ on $M$ is certainly integrable, and so we have by the
theorem a complete integral for the vector field $Z^e$ on the tangent
bundle. We describe the integral explicitly (this then also describes
the vector field, by restriction): it is the map
$R\times M^D \to M^D$ given by $(t, \beta ) \mapsto [d\mapsto \beta
(e^{-t}\cdot d)]$.---  The vector field on $M^D$ obtained this way
is, except for the sign, the {\em Liouville vector field}, cf.\ \cite
{G}, IX.2.

\section{Generalities on distributions}

We want to apply parts of the general theory of ordinary differential
equations to some of the basic equations of mathematical physics, the
wave- and heat- equations. This takes us by necessity to the realm of
distributions. Not  primarily as a technique, but because of the
nature of these equations: they model evolution through time of (say)
a {\em heat distribution}. A heat distribution is an {\em extensive}
quantity, and does not necessarily have a density {\em function},
which is an {\em intensive} quantity; the most important of all
distributions, the point distributions (or Dirac distributions), for
instance, do not. For the case of the heat equation, it is well known
that the evolution through time of any distribution ``instantaneously"
(i.e., after any {\em positive} lapse of time, $t>0$) leads to
distributions that do have smooth density functions. But in SDG, we
are interested also in what happens after a nilpotent lapse of time.
In more computational terms, we are interested in the Taylor
expansion of the solutions of evolution equations. For this, it is
necessary to stay within one vector space, that of distributions.

The vector space of ``distributions of compact support" on any object
$M$
can be introduced purely synthetically (see \cite{qr} p. 393, or
\cite{mr} p. 94)  as the $R$-linear dual of the vector space $R^M$
(which internally represents the vector space of smooth functions on
$M$). What follows could, to a certain extent (in particular for the
wave equation), be treated purely synthetically.

Presently, we shall only be interested in distributions on $R, R^2$
, and $R^3$, so for the presentation, we have chosen to assume that
we are working in a sufficiently good ``well-adapted" model ${\bf E}$
of SDG, containing the category of smooth manifolds as a full
subcategory. In such models, for any given manifold $M$, we could define
the linear subspace ${\cal D}(M)$ of $R^M$  consisting of functions
with compact support, (the ``test functions"). Then the vector space
of distributions  on $M$, ${\cal D}'(M)$, is taken to be the
$R$-linear dual of ${\cal D}(M)$.

One could take an alternative, slightly more concrete, approach:
namely,  take a model ${\bf E}$ of SDG which contains  the category
of smooth manifolds as above, but which also contains the category of
Convenient Vector Spaces
\cite{fro} and the smooth maps between them as a full subcategory.
The embedding is to preserve the cartesian closed structure. Such
models do exist: we provided in \cite{k2}, \cite{kr1} such an
embedding of Convenient Vector Spaces into the ``Cahiers" topos of
Dubuc \cite{dub}. Note that the usual topological (Fr\'{e}chet) vector
spaces of smooth functions, test functions, distributions, etc. on a
smooth manifold $M$ have canonical structure of Convenient Vector Spaces.  In
such a model, we can construct internal functions, say curves $f: R\to
{\cal D}'(M)$, by constructing, externally, a function by an ``excluded
middle" recipe of the form $$f(t)= ..\mbox{ if } t\neq 0\mbox{ ; }
f(t)= ..  \mbox{ if }t =0,$$ and then proving smoothness of $f$ by a
usual limit argument.

We have to resort to this kind of ``external" constructions only for
the heat equation, and there our embedding from \cite{k2}, \cite{kr1}
is not quite good enough, since it does not take manifolds with
boundary into account; for the heat equation, one constructs
externally an ``evolution"  map
$$R_{\geq 0} \to {\cal D}'(R)$$
by an excluded middle recipe.

So, for the justification of our treatment of the heat equation, we
need an extension (hopefully forthcoming) of our work
\cite{k2}, \cite{kr1}, i.e., we need  to construct a Cahiers-like
topos that includes also manifolds {\em with} boundary, and then to
construct an  embedding of Convenient Vector Spaces into that ``extended" 
Cahiers Topos.  (Maybe even the Cahiers Topos itself will be good 
enough.)

For what follows about wave equation, the Cahiers Topos, and the
embedding of Convenient Vector Spaces into it, is sufficient; in
fact, for these equations, a purely synthetic treatment alluded to
will be sufficient, since the distributions considered there are all
of compact support.

As stressed by Lawvere in \cite{law}, distributions should {\em not}
be thought of as generalized functions: functions are intensive
quantities, and transform {\em contravariantly}; distributions are
extensive quantities and transform {\em covariantly}. For functions,
this is the fact that the ``space" of functions on $M$, $R^M$ is
contravariant in $M$, by elementary cartesian-closed category theory.
Similarly, the ``space" of  distributions of compact support on $M$
is a subspace of $R^{R^M}$ (carved out by  the $R$ linearity
condition), and so for similar elementary reasons is covariant in
$M$. We shall write ${\cal D}'_c (M)$ for this subspace. The space of
functions of compact support on $M$ is only functorial with respect
to {\em proper} smooth maps, (counterimages of compact set required
to be compact), and so similarly, the space ${\cal D}' (M)$ of all
distributions on $M$ is covariant functorial only w.r.to proper maps.
The formula for covariant functorality looks the same  for ${\cal D}'$
and ${\cal D}' _c$; let us make it explicit for the ${\cal D}'$ case.
Let $f:M\to N$ be a proper map. The map 
${\cal D}'(f):{\cal D}'(M) \to {\cal D}'(N)$
is described by declaring
\begin{equation}<{\cal D}'(f)(\mu ) , \phi > = <\mu , \phi \circ f >,
\label{functorality}\end{equation}
where   $\mu$ is a distribution on $M$, and $\phi$ is a test function
on $N$, (so $\phi \circ f$ is a test function on $M$, by properness
of $f$). The brackets denote evaluation of distributions on test
functions.

We shall also write just $f(\mu )$ instead of ${\cal D}'(f)(\mu )$.

\medskip

Recall that a distribution $\mu$ on $M$ may be mulitplied by any
function $g:M\to R$, by the recipe
\begin{equation}<g\cdot \mu , \phi > = <\mu , g\cdot \phi >,
\label{multiplying}\end{equation}
observing that $g\cdot \phi$ is a test function (has compact support)
if $\phi$ is.

If $X$ is a vector field on $M$, one defines the directional
derivative $D_X (\mu )$ of a distribution $\mu$ on $M$ by the formula
\begin{equation}<D_X (\mu ), \phi > = -<\mu , D_X (\phi )>.
\label{Lie-d}\end{equation}
This in particular applies to the vector field
$\;{\partial}/{\partial x}\;$ on $R$, and reads here

\noindent $<\mu ',\phi > = -<\mu , \phi ' >$ ($\phi '$ denoting the
ordinary derivative of the function $\phi$). One has the following
Leibniz rule:

\begin{equation}
D_X (f\cdot \mu )= D_X (f) \cdot \mu + f\cdot D_X (f)
\label{Leibniz}\end{equation}
for any distribution $\mu$ and function $f $ on $M$. This is an
elementary consequence of the Leibniz rule for directional
derivatives $D_X$ of functions on $M$.

\medskip

\noindent {\bf Remark.} The equation (\ref{Lie-d}) becomes a theorem,
rather than a definition, if one takes the following line of
reasoning:
let $F$ be a covariant functor from microlinear spaces (and
invertible maps  between them)  to Euclidean vector spaces.
 Then one may define the {\em Lie derivative} along $X$, $L_X (\alpha
)$, as a map $F(M)\to F(M)$.  For the functor $F={\cal D}'$, $L_X$
becomes the $D_X$ described. We shall not pursue this line further
here.

\medskip
Applying $D_X$ twice leads to
$$<D_X (D_X (\mu )), \phi > = < \mu , D_X (D_X (\phi ))>$$
In particular, for $\mu$ a distribution on $R^n$
$$<{\partial ^2}/{\partial x_i}^2 (\mu) , \phi > = <\mu , 
{\partial ^2}/{\partial x_i}^2 (\phi )>$$ and therefore for the Laplace
operator
$\Delta = \sum \partial ^2 / \partial^2 x_i = \mbox{ div} \circ \mbox{
grad}$, we put
\begin{equation}<\Delta (\mu ),\phi > = <\mu , \Delta (\phi
)>.\label{Laplace}\end{equation}

\medskip
The following Proposition is an application of the covariant functorality of
the functor ${\cal D}_c$, which  will be used in connection with the
wave equation in dimension 2. We consider the (orthogonal) projection
$p: R^3 \to R^2$ onto the $xy$-plane. (It is not a proper map, so
functorality only works for compactly supported distributions.)

\begin{proposition} For any distribution $S$ (of compact support) on
$R^3$,
$$p(\Delta (S)) = \Delta (p(S)).$$
\label{proj}\end{proposition}
(The same result holds for any orthogonal projection $p$ of $R^n$ onto
any linear subspace; the proof is virtually the same, if one uses
invariance of $\Delta$ under orthogonal transformations.)

\medskip

{\bf Proof.} Let $\psi$ be any test function on $R^2$. Then
\begin{equation}<p(\Delta (S)),\psi > = < \Delta (S), \psi \circ p>
=  <S, \Delta (\psi \circ p)>\label{lap1}\end{equation}
But, with $\psi = \psi (x,y)$, $\psi \circ p$ is just $\psi$,
considered as a function of $x,y,z$ which happens not to depend on
$z$; so
$$\Delta (\psi \circ p )  =\frac{\partial \psi}{\partial x} +
\frac{\partial \psi}{\partial y} + \frac{\partial \psi}{\partial z};$$
the last term vanishes because $\psi$ does not depend on $z$, so the
equation continues
$$\;\;\;\;\;\;\;\;\;\;\;\;\;\;\;\;\;\;\;\;\;\;\;\;\;= \frac{\partial \psi}{\partial x} + \frac{\partial \psi}{\partial
y} = (\Delta (\psi ))\circ p.$$ So the right hand expression in
(\ref{lap1}) may be rewritten as
$$= <S, \Delta (\psi )\circ p > = <p(S), \Delta (\psi )> = <\Delta
(p(S)), \psi>,$$ from which the result follows.

\subsection{Spheres and balls as distributions}
For $a, b \in R$, we let $[a,b]$ denote the distribution
$f\mapsto \int _a ^b f(x)\; dx$. Such distributions on the line, we
of course call {\em intervals}; the {\em length} of an interval
$[a,b]$ is defined to be $b-a$. Note that the interval $[a,b]$ as a
distribution is not quite the same as the order theoretic interval,
i.e., the subset of $R$ consisting of $x$ with $a\leq x \leq b$. For
instance, the order theoretic interval from $0$ to $0$ contains all
nilpotent elements, whereas the distribution $[0,0]$ is the zero
distribution. The distribution theoretic interval $[a,b]$ contains
more information about $a$ and $b$ than does the order theoretic one.
We consider the question to which extent $[a,b]$ determines the
endpoints. The answer is contained in

\begin{proposition}Let $[a_1 ,b_1]$ and $[a_2 ,b_2]$ be two intervals 
in the distribution theoretic sense. They are equal as distributions 
if and only if they have same length, $b_1 - a_1 = b_2 - a_2$ ($=l$, say), and  $l\cdot (a_1 - a_2 ) = 0$ (this then also
implies $l\cdot (b_1 - b_2 )=0)$.\end{proposition}

{\bf Proof.} Assume $[a_1 ,b_1 ]=[a_2 ,b_2 ]$. The statement about
length follows immediately by applying each of these two
distributions to the function $f$ which is constant $1$. Generally, we
have for any function $f$ that
$$\int _{a_1} ^{b_1} f(x)\; dx =\int _{a_2} ^{b_2} f(x)\; dx$$
$$= \int_{a_1} ^{b_1} f(t+ a_2 - a_1)\; dt,$$
by making the change of variables $t= x+ a_1 - a_2 $. Subtracting, we
get
$$0= \int _{a_1} ^{b_1} (f(x) - f(x + a_2 - a_1 ))\; dx.$$
Apply this equation to the function $f(x)=x$, we get
$$0=\int _{a_1} ^{b_1}  (x-(x + a_2 - a_1 ))= (a_1 - a_2)\cdot (b_1 -
a_1) =
(a_1 - a_2)\cdot l.$$
Conversely, assume  $b_1 - a_1 = b_2 - a_2$ ($=l$, say), and $0= l\cdot (a_1
- a_2)$. For any function $f$, we calculate the values of the
distribution $\;[a_1 ,b_1 ]$ on $f$. We  have
$$  [a_1 ,b_1 ]  (f) = (b_1 - a_1 ) \int _0 ^1
f(a_1 + t\cdot (b_1 - a_1 ))\; dt= l  \int _0 ^1 f(a_1 + t\cdot l)\; dt.$$

Similarly
$$[a_2 ,b_2 ]  (f) =l \int _0 ^1 f(a_2 + t\cdot l)\; dt.$$
The difference is
\begin{equation}l\int _0 ^1 (f(a_1 + t\cdot l) - f(a_2 + t\cdot l ))\;
dt.\label{Had}\end{equation}
By Hadamard's Lemma, $f(a_1 + t\cdot l) - f(a_2 + t\cdot l )$ may be
written as
$(a_1 - a_2 )\cdot g(a_1 , a_2 , t)$ for some function $g$, and so
the integral (\ref{Had}) can be written as

$$= l\cdot (a_1 - a_2) \int _0 ^1 g(a_1 , a_2 ,t )\; dt,$$
which vanishes if $l\cdot (a_1 - a_2) =0$.

The assertions about $b_1 - b_2$ is similar.

\medskip

Note the following Corollaries: First, if the length $b_1 - a_1 $ of
an interval $[a_1 ,b_1 ]$ is invertible (positive, say), then the
endpoints $a_1$, $b_1$ are uniquely determined by the distribution $[a_1
,b_1 ]$.
Secondly, for any $t_1 , t_2$, we have
$$[-t_1 , t_1 ]= [-t_2 , t_2 ] \mbox{ implies } t_1 = t_2 .$$
In fact, by the Proposition, their lengths must be equal, i.e., $2t_1 = 2
t_2$.
The distribution $[-t ,t ]$ will appear below under the name $B_t$, ``the
ball of radius $t$ in dimension One".

\medskip
We shall also consider such ``balls" in dimension Two and Three,
where, however, $t$ cannot in general be recovered from the
distribution, unless $t$ is strictly positive.

We fix a positive  integer $n$. We shall
consider the sphere $S_t$ of radius $t$, and the ball $B_t$ of radius
$t$, for any $t\in R$, as distributions on $R^n$ (of compact support,
in fact), in the following  sense:
$$<S_t , \psi > = \int _{S_t} \psi (x) dx = t^{n-1} \int _{S_1} \psi
(t\cdot u)\; du,$$

$$<B_t , \psi > = \int _{B_t} \psi (x) dx = t^{n} \int _{B_1} \psi
(t\cdot u)\; du,$$
where $du$ refers to the surface element of the unit sphere $S_1$ in
the first equation and to the volume element of the unit ball $B_1$ in
the second. The expressions involving $\int _{S_t}$ and $\int_{B_t}$
are to be understood symbolically, unless $t>0$; if $t>0$, they make
sense literally as integrals over sphere and ball, respectively, of
radius $t$, with $dx$ denoting surface-, resp.\ volume element. But the
expression on the right in both equations make sense for any $t$, and
so the distributions $S_t$ and $B_t$  are defined for all $t$; in
particular, for nilpotent ones.

It is natural to consider also the following distributions $S^t$ and
$B^t$ on $R^n$ (likewise of compact support):
$$<S^t , \psi > = \int _{S_1} \psi (t\cdot u)\; du,$$
$$ <B^t , \psi > = \int _{B_1} \psi (t\cdot u)\; du.$$
For $t>0$, they may, modulo  factors of the type $4\pi$,  be
considered as ``average over $S_t$" and ``average over $B_t$",
respectively, since $S^t$ differs from $S_t$ by a factor $t^{n-1}$,
which is just the surface area of $S_t$ (modulo the factor of type
$4\pi$), and similarly for $B^t$.

Note that $S^1 = S_1$ and $B^1 = B_1$. And also note that the
definition of $S^t$ and $B^t$ can be formulated as
$$S^t = H_t (S_1 )\mbox{ , } B^t = H_t (B_1),$$
where $H_t : R^n \to R^n$ is the homothetic transformation $u\mapsto
t\cdot u$, and where we are using the covariant functorality of
distributions of compact support.

\medskip

For low dimensions, we shall describe the distributions $S_t$, $B_t$, $S^t$
and $B^t$  explicitly:

\medskip
\noindent{\bf Dimension 1}
$$<S_t , \psi > = \psi (-t) + \psi (t)$$
$$<B_t , \psi > = \int _{-t} ^t \psi (s)\; ds$$
$$<S^t , \psi > = \psi (-t) + \psi (t)$$
$$<B^t , \psi > = \int _{-1} ^1 \psi (t\cdot s)\; ds$$

\medskip
\noindent{\bf Dimension 2}
$$<S_t , \psi > =\int _0 ^{2\pi} \psi (t\cos \theta , t\sin \theta
)\; t\; d\theta$$
$$<B_t , \psi > =\int _0 ^t \int _0 ^{2\pi} \psi (s \cos \theta ,
s\sin \theta) \; s\; d\theta \; ds $$
$$<S^t , \psi > =\int _0 ^{2\pi} \psi (t\cos \theta , t\sin \theta
)\; d\theta $$
$$<B^t , \psi > =\int _0 ^1 \int _0 ^{2\pi}\psi (t s\cos \theta , t\
s\sin \theta )\; s\; d\theta \; ds $$

\medskip
\noindent{\bf Dimension 3}
$$<S_t , \psi > =\int _0 ^{\pi} \int _0 ^{2\pi}\psi (t\cos \theta
\sin \phi , t\sin \theta \sin \phi , t\cos \phi ) t^2 \sin \phi \;
d\theta \; d\phi $$
$$<B_t , \psi > =\int _0 ^t \int _0 ^{\pi} \int _0 ^{2\pi}\psi (s\cos
\theta \sin \phi , s\sin \theta \sin \phi , s\cos \phi )\; s^2 \sin
\phi \; d\theta \; d\phi \; ds$$
$$<S^t , \psi > =\int _0 ^{\pi}\int _0 ^{2\pi} \psi (t\cos \theta
\sin \phi , t\sin \theta \sin \phi , t\cos \phi ) \sin \phi \;
d\theta \; d\phi $$
$$<B^t , \psi > =\int _0 ^1 \int _0 ^{\pi} \int _0 ^{2\pi}\psi
(ts\cos \theta \sin \phi , ts\sin \theta \sin \phi , ts\cos \phi )\;
s^2 \sin \phi \;d\theta \; d\phi \; ds.$$

Notice that these formulas make sense for all $t$ (positive,
negative, nilpotent, ...  ), using the standard convention :$\int _a ^b = - 
\int _b ^a$), whereas set-theoreti\-cally $S_t$ and $B_t$ (as point 
sets) only make good sense for $t>0$.

It is clear from the very definition that $S_t = t^{n-1} S^t$ 
and $B_t = t^n B^t$ (in any dimension $n$); but since we 
are interested also 
in $t$'s that are not invertible, $S_t$ and $S^t$ cannot be defined in
terms of each other.

Note also that $S_0 = B_0 = 0$, whereas $S^0$ and $B^0$ are constants 
 times the Dirac distribution at the origin $0$.  The 
constants are the ``area'' of the unit sphere, or the ``volume'' of the unit 
ball, in the appropriate dimension.  Explicitly,
$$S^0 =2\cdot \delta (0),\; 2\pi \cdot \delta (0),\; 4\pi \cdot \delta (0),$$
and
 $$B^0 = 2\cdot \delta (0),\; \pi \cdot \delta (0),\; 
\frac{4\pi}{3}\cdot \delta (0)$$
in dimensions 1,2, and 3, respectively.

 \medskip We 
shall also have occasion to consider the distribution (of compact
 support) $t\cdot S^t$ on $R^3$ as well as its projection $p(t\cdot 
S^t)$ on the $xy$-plane (using functorality of ${\cal D}'_c$ with 
respect to the projection map $p:R^3 \to R^2 $).  For $t>0$ (more 
generally, for $t$ invertible), we can give an explicit integral 
expression for it, but note that since $S^t$ and $t\cdot S^t$ are 
defined for all $t$, then so is $p(t\cdot S^t)$, whether or not we 
have such an integral expression.  The integral expression (for $t>0$) 
goes under the name of {\it Poisson kernel} for the wave equation in 
dimension 2 and may be obtained as follows: using the above expression 
for $S^t$ in dimension 3, we have for a test function $\psi$ that only 
depends on $x,y$, but not on $z$ that $$<t\cdot S^t , \psi > =\int _0 
^{\pi}\int _0 ^{2\pi} \psi (t\cos \theta \sin \phi , t\sin \theta \sin 
\phi )\cdot t \cdot \sin \phi
\; d\theta \; d\phi .$$
We then make the change of variables $\rho =t\sin \phi$,  $\phi =
\arccos \rho /t$, $d\phi = d\rho / \sqrt{t^2 - \rho ^2}$, and then
the integral becomes
$$2\int _0 ^{\frac{\pi}{2}} \int _0 ^ {2\pi} \frac{\psi (\rho \cos
\theta , \rho \sin \theta )\; \rho \; d\theta \; d\rho }{\sqrt{t^2 -
\rho ^2}},$$
using the explicit form of the ball distribution $B_t$ in dimension
2, we may rewrite the right hand side here as
$$<\frac{2}{\sqrt{t^2 - \rho ^2}}\cdot B_t , \psi >,$$
so that we have, for $t>0$ (or even for $t$ invertible),
\begin{equation}p(t\cdot S^t ) =\frac{2}{\sqrt{t^2 - \rho ^2}} \cdot
B_t .\label{poisson}\end{equation}

\section{Vector Calculus}
The Main Theorem of vector calculus is Stokes' Theorem:
$\;\int _{\partial \gamma} \omega = \int _{\gamma} d\omega$, for
$\omega$ an $(n-1)$-form, $\gamma$ a suitable $n$-dimensional
figure (with appropriate measure on it) and $\partial \gamma$ 
its geometric boundary. In the synthetic context, the theorem 
holds at least for any singular cubical chain  $\gamma :I^n 
\to M$ ($I^n$ the $n$-dimensional coordinate cube), because the 
theorem may then be reduced to the fundamental theorem of calculus, 
which is the only way integration enters in the elementary synthetic 
context; measure theory not being available therein.  For an account 
of Stokes' Theorem in this context, see \cite{mr} p.139.  Below, we 
shall apply the result not only for singular {\em cubes}, but also for 
singular {\em boxes}, like the usual $\gamma : [0, 2\pi ] \times [0,1] 
\to R^2$, parametrizing the unit disk by polar coordinates, 
\begin{equation}\gamma (\theta , r ) = (r\cos \theta , r \sin \theta 
).
\label{polar}\end{equation}
We shall need from vector calculus the Gauss-Ostrogradsky 
``Divergence Theorem"
$$\mbox{flux of } {\bf F} \mbox{ over } \partial \gamma = \int
_{\gamma} (\mbox{divergence of } {\bf F}),$$
with ${\bf F}$ a vector field, for the geometric ``figure" 
$\gamma $ = the unit ball in $R^n.$
 For the case of the unit ball in $R^n$, the reduction of the Divergence Theorem to Stokes' Theorem is  a matter of the {\em differential} calculus of vector fields, differential forms, inner products etc.\, (See e.g.  \cite{lang} p. 204). For the convenience of the reader, we recall the case $n=2$.

\medskip
Given a vector field ${\bf F}(x,y) = (F(x,y),\; G(x,y))$ in $R^2$, apply
Stokes' Theorem to the differential form
$$\omega := -G(x,y)dx + F(x,y)dy$$ 
for the singular rectangle $\gamma$ given by (\ref{polar}) above. Then 
$$
\left \{
\begin{array}{lll}
\gamma ^* (dx)= \cos \theta dr - r \sin \theta d\theta \\
\gamma ^* (dy ) = \sin \theta dr + r \cos \theta d\theta \\
\gamma^* (dx\wedge dy) = r \; (dr \wedge d\theta) 
\end{array}
\right.$$
Since $\;d\omega = (\partial G/\partial y + \partial F /
\partial x )\; dx\wedge dy = \mbox{ div }({\bf F}) \; dx\wedge dy$, then
$$\gamma^* (d\omega ) =\mbox{ div }({\bf F}) \; r \; (dr\wedge d\theta)
$$ 
On the other hand,
\begin{equation}\gamma ^* \omega = (F \; \sin \theta - G \; \cos
\theta )dr + (F \; r \; \cos \theta + G \; r \; \sin \theta )\;
d\theta ,\label{oneform}\end{equation}
(all $F$, $G$, and ${\bf F}$ to be evaluated ar $(r\cos \theta ,
r\sin \theta )$).
Therefore $$\int _{\gamma} d\omega = \int _0 ^{2 \pi}\int _0 ^1 \mbox{ div
}({\bf F}) \; r \; dr\;  d\theta ;$$
this is $\int _{B_1} \mbox{ div }({\bf F}) \;dA$.
On the other hand  by Stokes' Theorem
$\int _{\gamma} d\omega = \int _{\partial \gamma} \omega$ which is a
curve integral of the 1-form (\ref{oneform}) around the boundary of
the rectangle $[0, 2\pi ] \times [0,1]$. This curve integral is a sum
of four terms corresponding to the four sides  of the rectangle. Two
of these (corresponding to the sides $\theta = 0$ and $\theta = 2\pi
$) cancel, and the term corresponding to the side where $r=0$
vanishes because of the $r$ in $r\; (dr\wedge d\theta)$, so only the
side with $r=1$, $0\leq \theta \leq 2\pi$ remains, and its
contribution is, with the correct orientation,
$$\int _0 ^{2\pi} (F (\cos \theta , \sin \theta )\cos \theta + G (\cos
\theta , \sin \theta ) \sin \theta )\; d\theta = \int _{S_1} {\bf
F}\cdot {\bf n} \; ds$$ where ${\bf n}$ is the outward unit normal of
the unit circle. This expression is the flux of ${\bf F}$ over the
unit circle, which thus equals the divergence integral calculated
above.

\medskip

We insert for reference two obvious ``change of variables" equations.
Recall that $H_t : R^n \to R^n$ is the homothetic transformation
``multiplying by $t$". We have, for any vector field ${\bf F}$ on
$R^n$ (viewed, via principal part, as a map $R^n \to R^n $):
\begin{equation}\mbox{div } ({\bf F}\circ H_t ) = t\cdot (\mbox{div
}{\bf F})\circ H_t ,
\label{littlediv}\end{equation}
and
\begin{equation}t^n \int_{B_1} \phi \circ H_t = \int _{B_t} \phi .
\label{ch-var}\end{equation}

We now combine vector calculus with the calculus of the basic
ball- and sphere-distributions, as introduced in Section 3, to prove

\begin{theorem} In $R^n$ (for any $n$), we have, for any $t$,
$$\frac{d}{dt} S^t = t\cdot \Delta(B^t),$$
($\Delta =$ the Laplace operator).
\label{prop1}\end{theorem}

{\bf Proof.} We prove first that
$$t^{n-1}\cdot \frac{d}{dt} S^t = t^n \cdot \Delta(B^t) .$$
In fact, for any test function $\psi$,
$$< t^{n-1}\cdot \frac{d}{dt} S^t , \psi > = t^{n-1}\cdot
\frac{d}{dt} \int _{S_1} \psi (tu)\; du = t^{n-1} \int _{S_1} (\nabla
\psi) (tu)\; \cdot u \; du,$$
(by differentiating under the integral sign and using the chain rule)
$$= t^{n-1} \cdot \mbox{ flux of }((\nabla \psi)\circ H_t )\mbox{
over } S_1 ,$$
where $H_t :R^n \to R^n$ is the homothetic transformation
``multiplying by $t$". This, by the Divergence Theorem, may be
rewritten as

$$t^{n-1} \int_{B_1}  \mbox{ div } ((\nabla \psi )\circ H_t )  =
t^{n} \int_{B_1} (\mbox{ div }(\nabla \psi ))\circ H_t,$$
(using (\ref{littlediv}))
$$=t^n \cdot \int _{B_1} (\Delta \psi ) \circ H_t = \int _{B_t}
\Delta \psi $$
(by a standard change of variables, cf.\ (\ref{ch-var})), so
$$= <B_t , \Delta \psi > = < t^n \cdot B^t ,\Delta \psi > = t^n \cdot
< \Delta B^t , \psi >.$$

From
$$t^{n-1} < \frac{d}{dt} S^t , \psi > = t^n <\Delta B^t , \psi >,$$
we may of course conclude the desired equality, by cancelling
$t^{n-1}$ on both sides, if $t$ is invertible; but we want the
equation {\em for all} $t$. We can get this from ``Lavendhomme's
principle", which says that if $f:R\to R$ satisfies $t\cdot f(t)= 0$ 
for all $t$, then $f(t)$ is constantly $0$.  This principle was 
derived from the integration axiom purely synthetically by Lavendhomme 
in \cite{lav} p.25.  So the claim of the Theorem is valid for all $t$.

\medskip
We collect information about $t$-derivatives of the four basic
distributions $S_t$ , $B_t$, $S^t$ and $B^t$ in $R^n$. The results
are valid for any $n$ and any $t$. For invertible $t$ (say positive
$t$), some of the statements may be simplified by multiplying by
$t^{-1}$, but we prefer  having  formulae which are universally valid.

\begin{theorem} We have in dimension $n$ for all $t$:

\begin{equation}
\frac{d}{dt} (B_t ) = S_t ,
\label{a1}\end{equation}

\begin{equation}
t\cdot \frac{d}{dt}(S_t ) = (n-1)\;  S_t + t\cdot \Delta (B_t) ,
\label{a2}\end{equation}

\begin{equation}
t\cdot \frac{d}{dt} (B^t ) = S^t - n\;  B^t ,
\label{a3}\end{equation}

\begin{equation}
\frac{d}{dt}(S^t )= t\cdot \Delta (B^t ),
\label{a4}\end{equation}

In dimension $1$, we also have
\begin{equation}
\frac{d}{dt}(S_t ) = \Delta B_t
\label{a5}\end{equation}

\label{allfour}\end{theorem}

{\bf Proof.} Equation (\ref{a1}) is an immediate consequence of the
fundamental theorem of calculus; e.g.\ for $n=2$, consider the
explicit formula for $B_t$ given above in Section 3 (``Spheres and balls as 
distributions").  With $\int _0 ^t$ as the outer integral, the $d/dt$ 
of it is just the inner integral, i.e., exactly the exhibited formula 
(idem) for $S_t$.

For (\ref{a2}), we $t$-differentiate the equation $S_t = t^{n-1}\cdot
S^t $ by the Leibniz rule and get $(n-1)\cdot t^{n-2} \cdot S^t +
t^{n-1}\cdot d/dt (S^t )$; so by Theorem \ref{prop1}, 

$$d/dt (S_t ) = (n-1)\cdot t^{n-2} \cdot S^t +
t^{n-1}\cdot t\cdot \Delta (B^t ).$$ If we multiply this equation by
$t$, we get
$$t\cdot d/dt (S_t ) = (n-1)\cdot t^{n-1} \cdot S^t + t^{n+1} \cdot
\Delta (B^t );$$
using $S_t = t^{n-1} S^t$ and $B_t = t^n B^t$, the result follows
(note that $\Delta$ commutes with multiplication by $t$).

The proof of (\ref{a3}) is similar: $t$-differentiating  $  t^n
\cdot B^t =B_t$, we get
$$n\cdot t^{n-1} \cdot B^t + t^n \cdot d/dt B^t = d/dt B_t = S_t ,$$
(using (\ref{a1})), so using $S_t = t^{n-1} S^t$, this equation may
be rewritten as
$$t^{n-1}\cdot (n\cdot B^t  + t\cdot  d/dt B^t ) = t^{n-1}\cdot S^t
$$
The result now follows by cancelletion of the factor $t^{n-1}$ by
Lavendhomme's principle, and rearranging.

Next, (\ref{a4}) is identical to Theorem \ref{prop1}, and is included
again for completeness' sake.

Finally, (\ref{a5}) follows from (\ref{a2}): the first term vanishes,
since $n-1 =0$, and in the remaining equation, we may cancel the
factor $t$ by Lavendhomme's principle. Alternatively, (\ref{a5}) can
be proved directly, by a very simple calculation.

\section{Wave equation}
Let $\Delta$ denote the Laplace operator $\sum \partial ^2 /\partial 
x_i ^2$ on $R^n$. We shall consider the wave equation (WE) in $R^n$, 
(for $n=1,2,3$), 
\begin{equation}\frac{\partial ^2}{\partial t^2} Q = \Delta Q
\label{waveequation}\end{equation}
as a second order ordinary differential equation on the Euclidean
vector space ${\cal D}' _c (R^n )$ of distributions of compact support;
in other words, we are looking for functions $$Q: R \to {\cal D}' _c (R^n
)$$ so that for all $t\in R,$ $\ddot{Q}(t) = \Delta (Q(t))$ (viewing
$\Delta $ as a map $ {\cal D}' _c (R^n ) \to {\cal D}' _c (R^n ) $.  We
shall only be looking for {\em particular } solutions, in fact, so
called {\em fundamental} solutions: solutions whose initial value and 
initial speed is either the Dirac distribution at $0$, or $0$.  Given any other 
initial value and speed --- these being both assumed to be 
distributions of compact support ---, the corresponding particular 
solution may, as is well known, be obtained from the fundamental 
solution just by {\em convolution} $*$ with these fundamental solutions.  This 
follows purely formally from the rules for convolution of 
distributions $P$ and $Q$, such as $Q * \delta (0) = Q$, $D(P * Q )= 
D(P ) * Q$, where $D$ is any differential operator on $R^n$ with 
constant coefficients; and from linearity of the convolution, implying 
that $d/dt (P_t * Q )= (d/dt\; P_t )* Q$; see e.g.  \cite{sch}, Ch.  
3.  
\medskip

\noindent{\bf Dimension 1}

\begin{theorem} The function $R \to {\cal D}' _c (R)$ given by $$t\mapsto
1/2 \cdot S^t $$ is a solution of the WE in dimension 1; its initial
value and speed are, respectively $\delta (0)$ and $0$.

The function $R \to {\cal D}' _c (R)$ given by $$t\mapsto 1/2 B_t$$ is a
solution of the WE; its initial value and speed are, respectively, $0$
and $\delta (0)$.
\label{dim1}\end{theorem}

{\bf Proof.} The statements about the initial values are immediate
from the explicit integral formulas for $B_t$ and $S^t$ (putting
$t=0$). The statements about the initial speeds are equally immediate
from the following formulas (\ref{b1}) and (\ref{b3}) for the
$t$-derivatives, (putting $t=0$). We have by (\ref{a4})

\begin{equation}\frac{d}{dt} (S^t ) = t\cdot \Delta (B^t
),\label{b1}\end{equation}
and so by further $t$ differentiation
$$\frac{d^2}{dt^2} (S^t ) = \Delta (B^t ) + t\cdot \frac{d}{dt}
(\Delta (B^t ));$$
now, $d/dt$ and $\Delta $ commute, so we may continue
$$= \Delta (B^t ) + \Delta (t\cdot \frac{d}{dt}  B^t ) = \Delta (B^t
) + \Delta (S^t - 1\cdot B^t ),$$
using (\ref{a3}) with $n=1$. Now by linearity of $\Delta$, the terms
involving $B^t$ in the last expression cancel, and we are left with

\begin{equation}\frac{d^2}{dt^2} (S^t ) = \Delta (S^t ),
\label{b2}\end{equation}
which establishes WE for $S^t$ and hence also for $1/2 \cdot S^t$.

Also, by (\ref{a1}), we have that
\begin{equation}\frac{d}{dt} (B_t ) = S_t ,\label{b3}\end{equation}
and so by further $t$ differentiation
$$\frac{d^2}{dt^2} (B_t ) = \frac{d}{dt} (S_t ) = \Delta (B_t ),$$
using (\ref{a5}), which establishes WE for $B_t$ and hence for $1/2
\cdot B_t$. So the theorem is proved.

\medskip
\noindent {\bf Dimension 3} 
\begin{theorem} The function $R \to {\cal D}' _c 
(R^3 )$ given by $$t\mapsto \frac{1}{4\pi} \cdot t \cdot S^t $$ is a 
solution of the WE in dimension 3; its initial value and speed are, 
respectively, $0$ and $\delta (0)$.

The function $R \to {\cal D}' _c (R^3 )$ given by $$t\mapsto  \frac{1}{4\pi} \cdot (S^t
+ t^2 \cdot \Delta (B^t ))$$ is a solution of the WE; its initial value
and speed are, respectively, $\delta (0)$ and $0$.

\label{dim3}\end{theorem}

{\bf Proof.} We calculate first $d/dt$ of $t\cdot S^t$, using
(\ref{a4}):
\begin{equation}\frac{d}{dt}( t\cdot S^t
)= S^t + t^2 \cdot \Delta (B^t ),\label{b33}\end{equation}
and so by Theorem \ref{prop1} (= (\ref{a4})), $$\frac{d^2 }{dt^2 } 
(t\cdot S^t ) =t\cdot \Delta (B^t ) +2 \cdot t
\cdot \Delta (B^t ) +t^2 \cdot \Delta (\frac{d}{dt} B^t )$$
$$= 3\cdot t \cdot \Delta (B^t ) + t\cdot \Delta (t\cdot \frac{d}{dt}
B^t )$$
$$\;\;\;= 3\cdot t \cdot \Delta (B^t ) + t\cdot \Delta (S^t - 3B^t ),$$
using (\ref{a3}), and now by linearity of $\Delta$, the terms
involving $\Delta (B^t )$ cancel, so we are left with the equation

\begin{equation}
\frac{d^2 }{dt^2 } (t\cdot S^t )=\Delta (t\cdot S^t ),
\label{b4}\end{equation}
which establishes WE for $t\cdot S^t$ and hence for $1/4\pi \; \cdot t\cdot
S^t$.  The statements about initial value and speed are immediate 
(using (\ref{b33}) for the speed).

Because $d^2/dt^2$ and $\Delta$ commute, it is clear that if
$t\mapsto Q(t)$ is a distributional solution of WE, then so is
$t\mapsto d/dt Q(t)$. So since $t\cdot S^t$ is a solution, then so is
its $t$-derivative (calculated in (\ref{b33}) above), i.e.\
$S^t + t^2 \cdot \Delta (B^t )$ is a solution. Its initial value 
 and its initial speed can be found by
putting $t=0$ in (\ref{b4}) (note $\delta$ commutes
with multiplication by $t$).

\medskip

\noindent{\bf Dimension 2}

\medskip
Recall that we considered the orthogonal projection $p: R^3 \to R^2.$
Applying covariant functorality, we get for any distribution $Q$ on
$R^3$ of compact support a distribution $p(Q)$ on $R^2$, also of compact
support.

\begin{theorem} The function $R \to {\cal D}' _c (R^2 )$ given by
$$t\mapsto  \frac{1}{4\pi} \cdot p(t \cdot S^t )$$ 
is a solution of the WE in dimension
2; its initial value and speed are, respectively, $0$ and $\delta
(0)$.

The function $R \to {\cal D}' _c (R^2 )$ given by $t\mapsto 1/4\pi  
\cdot p(S^t + t^2 \cdot \Delta (B^t ))$ is also a solution of the WE 
in dimension 2; its initial value and speed are, respectively, $\delta 
(0)$ and $0$.
\label{dim2}\end{theorem}

Recall that an explicit integral formula for $p(t\cdot S^t )$, for
$t>0$, was given above, in (\ref{poisson}) (``Poisson kernel").

\medskip

{\bf Proof.} The fact that the two distributions in question are
solutions of the WE is  immediate from the Proposition \ref{proj}
(``$p$ commutes with $\Delta $") and from the  fact that ${\cal D}' _c
(p): {\cal D}' _c (R^3 )\to {\cal D}' _c (R^2 )$ is linear, and hence
commutes with formation of $d/dt$; also,  ${\cal D}' _c (p)$ sends 
Dirac distribution at $0\in R^3$ to Dirac distribution at
 $0\in R^2$, so the initial values and speeds are as claimed.

\medskip
The Taylor Series at $t=0$ for the solutions given can be calculated
directly, but they can more easily be obtained from the formal
solution given in Proposition \ref{exp-sol2}.

\section{Heat equation.}
In this section we deal with distributions that do not have compact
support and we only consider the one-dimensional case. We are thus
considering solutions for the vector field on the Euclidean vector
space ${\cal D}' (R)$, whose principal part is given by $\Delta :{\cal
D}'(R)\to {\cal D}' (R)$. We consider the particular solution $K:
R_{\geq 0} \to {\cal D}' (R)$ whose initial value is the distribution
$\delta (0).$ Thus, referring to the general treatment of solutions
for (differential equations given by) vector fields, we are
considering $\tilde{R} = R_{\geq 0}$; for the heat equation, one
cannot do better, as is well known.  Also, as mentioned above, we
rely on external (classical) calculus; namely, we consider the
classical ``heat kernel'' function, i.e., the function $K: R_{\geq 0}
\to {\cal D}' (R)$ given by
\begin{equation} K(t) = \left \{ \begin{array}{ll}
\frac{1}{\sqrt{4\pi t}}\; e^{-\frac{x^2}{4t}}&\mbox{ for } t>0\\
& \\
\delta (0)&\mbox{ otherwise }.
\end{array}
\right. \label{heatequation}\end{equation}

Here, for the case $t>0$, we described a function rather than a
distribution, so here we do make the identification of functions
$g(x)$ with distributions 
$\phi \mapsto \int_{-\infty}^{\infty} g(x) \phi (x) \; dx$. 
Differentiation of distributions reduces to differentiation of the 
representing functions.  For $t>0$, we thus have $K_t (x) = K(t,x)$, a
smooth function in two variables, described by the above expression.
It  satisfies the heat equation
$$\frac{\partial K}{\partial t} = \frac{\partial ^2 K}{\partial
x^2},$$
for $t>0$. Also the following limit expression is classical:
\begin{equation}
\mbox{lim}\; _{t\rightarrow 0^+}\int _{-\infty} ^{\infty} K(t,x) \phi
(x) \;
dx = \phi (0)
\label{limit1}\end{equation}
for any test function $\phi$. More generally,

\begin{proposition} For any integer $n\geq 0$, and any test function
$\phi$
\begin{equation}
\mbox{{\em lim}}\; _{t\rightarrow 0^+}\;\frac{\partial ^n}{\partial t^n}\int
_{-\infty} ^{\infty} K(t,x) \phi (x) \; dx = \phi ^{(2n)} (0).
\label{limit2}\end{equation}
\label{limit3}\end{proposition}

{\bf Proof.} The case $n=0$ is just (\ref{limit1}); the general case
follows by iteration. Let us do the case $n=1$. Then
\begin{eqnarray*}
\frac{\partial }{\partial t} \int _{-\infty} ^{\infty} K(t,x) 
\phi (x)\; dx & = & \int _{-\infty} ^{\infty}\frac{\partial}{\partial t}
 K(t,x) \phi (x)\;dx  \\ 
 & = & \int _{-\infty} ^{\infty} \frac{\partial ^2}{\partial x^2}
K(t,x)\phi (x)\; dx 
\end{eqnarray*}
(by the heat equation for $K$)
\begin{eqnarray*}
& \;\;\;\;\;\;\;\;\;\;\;\;\;\;\;\;\;\;\;\;\;\;\;\;\;\;\;\;\;\;\;= 
& \int _{-\infty} ^{\infty} K(t,x)\phi ^{(2)} (x)\; dx \\ 
\end{eqnarray*}
(by integration by parts.)

\medskip
We then use (\ref{limit1}), for the test function $\phi ^{(2)}$ to conclude 
(\ref{limit2}) for $n=1$.

\begin{proposition} The function $K: R_{\geq 0} \to {\cal D}' (R)$  is
smooth.
\label{newsmooth}\end{proposition}

Here, smoothness is taken in the following sense (appropriate for
convenient vector spaces): for each test function $\phi$, the function $R_{\geq
0} \to R$ given by $t\mapsto <K(t),\phi >$ is smooth.

\medskip
\noindent{\bf Proof.} It suffices to prove that $K$ is infinitely
often differentiable at $0$, since smoothness for $t>0$ is clear.  For fixed
$t>0$, we let $K_t$ denote the function in $x$ described in (the first
clause in) (\ref{heatequation}) above. Thus,  $<K(t), \phi >$ is given by the integral
\begin{equation}\int _{-\infty}^{\infty} \frac{1}{\sqrt{4\pi t}}\;
e^{-\frac{x^2}{4t}} \cdot \phi (x)\; dx.\label{ok!}\end{equation}

\medskip
We first notice that, by Hadamard's Lemma, $\phi (x)=\phi (0) + x\psi (x)$.
By linearity, $<K_t, \phi>=<K_t,\phi (0)> + <K_t, x\psi (x)>.$
But $ <K_t,\phi (0)>=\phi (0)$ and this implies that the derivative
of $<K_t, \phi>$ at $0$ is
\begin{equation}
\lim _{t\to 0^+}(1/t)<K_t, x\psi (x)>
\label{lim}
\end{equation}
To compute this limit, we use the formulas and notations in Lang's
book \cite{lang}, with the exception that we use $\cal{F}$ for the
Fourier transform.  We also use the following well known formulae, where all
the functions under considerations belong to the class $S$ of fast
decreasing functions and thus $\cal{F}$ works with no limitations.
First, for any pair of functions $\alpha$, $\beta$ in this class, one
has the ``adjointness'' formula
$$\int_{-\infty}^{\infty} \cal{F}(\alpha ) \beta =
\int_{-\infty}^{\infty} \alpha \cal{F}(\beta ).$$
Furthermore

\begin{equation}{\cal
F}((1/\sqrt{4\pi})e^{-t{\xi}^2})(x)=(1/\sqrt{4\pi
t})e^{-x^2/4t}\label{f-one}\end{equation}
 \begin{equation}{\cal{F}}(x\psi
(x))(\xi)=i({\cal{F}}({\psi}))'(\xi)\label{f-two}\end{equation}
\begin{equation}\xi{\cal F}(\psi)(\xi)=-i{\cal
F}({\psi}')(\xi)\label{f-three}\end{equation}

To show the existence of the limit, we compute, using (\ref{f-one}),
adjointness, and (\ref{f-two})
\begin{eqnarray}
<K_t, x\psi (x)> & = &
\int_{-\infty}^{\infty}(1/\sqrt{4\pi t})e^{-x^2/4t}x\psi (x) dx
\nonumber \\
 & = & \int_{\infty}^{\infty}{\cal F}((1/\sqrt{4\pi})e^{-t{\xi}^2})(x)
x\psi (x) dx  \nonumber \\
  & = &
\int_{\infty}^{\infty}((1/\sqrt{4\pi})e^{-t{\xi}^2}){\cal{F}}(x\psi
(x))(\xi)d\xi   \nonumber \\
  & = &
\int_{\infty}^{\infty}((1/\sqrt{4\pi})e^{-t{\xi}^2})i({\cal{F}}({\psi}))'(\xi)
d
\xi  \nonumber \\
  & = &
2it\int_{-\infty}^{\infty}
( 1/\sqrt{4\pi}) e^{-t{\xi}^2}\xi{\cal{F}}(\psi)(\xi)d\xi ,
\nonumber \end{eqnarray} The last step uses integration by parts.
Using (\ref{f-three}), this may be rewritten as

$$\;\;=2t \int _{-\infty} ^{\infty}
( 1/\sqrt{4\pi }) e^{-t\xi ^2} {\cal F}(\psi ')(\xi )\; d\xi$$
$$=
2t \int _{-\infty} ^{\infty} ( 1/\sqrt{4\pi t}) e^{-x^2 /4t} \psi
'(x) \; dx$$
using adjointness and (\ref{f-one}) in the last step.  Now we divide by
$t$, as requested in (\ref{lim}), and let $t\to 0^+$. Using (\ref{limit1}),
  we thus get that the limit in (\ref{lim}) equals 
$$\lim_{t\to 0^+} 2 \int _{-\infty}
^{\infty} ( 1/\sqrt{4\pi t}) e^{-x^2 /4t } \psi ' (x) \; dx = 2\cdot \psi ' (0).
$$
But since $\phi (x)= \phi (0) + x\cdot \psi (x)$, $2\cdot \psi '
(0) = \phi '' (0)$, This proves that the limit in (\ref{lim}) exists
and equals $\phi '' (0)$; we conclude that \begin{equation} \lim _{t\to
0^+}(1/t) [<K_t, \phi > -\; \phi (0)] = \phi''(0).
\label{socalledlemma}\end{equation}

\medskip
To better understand what has been done and to develop this matter
further, let us define for every $t\geq 0$
$$ f(t)=<K_t,\phi>$$
We can summarize the results of this section as follows
$$f'(0)= {\phi}''(0)$$
Recall from Proposition \ref{limit3} that
$$f^{(n)}(t)=<K_t, {\phi}^{(2n)}>$$
and thus, by going to the limit when $t\to 0^+,$
$$\lim_{t\to 0^+}f^{(n)}(t)= {\phi}^{(2n)}(0)$$
These results suffice to summarize the result in the present Section
in the following way:

\begin{corollary}
 The function $f$ is smooth and, furthermore, 
$f^{(n)}(0)={\phi}^{(2n)}(0).$
\end{corollary}
{\bf Proof.} Let us show, for instance, that $f''(0)$ exists and
equals ${\phi}''''(0).$

\noindent Using the previous results,
$(1/t)[f'(t)-f'(0)]=(1/t)[<K_t,{\phi}''> -\; {\phi}''(0)]$
and this implies the corollary, by going to the limit when $t\to 0^+$
and using (\ref{socalledlemma}) with the function ${\phi}''$ instead
of $\phi.$ Now, iterate.

\medskip
The idea to use Fourier transform to prove smoothness was pointed out
to us by H. Stetk\ae r and E. Skibsted.

\medskip
Summarizing: we have a smooth function $K: R_{\geq 0} \to
 {\cal D}'(R)$, satisfying the heat equation $\partial K / \partial t (t) =
\Delta (K(t))$ for all $t\geq 0$; for $t=0$, this follows from
Proposition \ref{limit3}.  By the assumed fullness of the embedding of
smooth manifolds with boundary and convenient vector spaces into the
model of SDG, we have the desired solution internally in the
model. We may then ask for the values of $K$ for {\em nilpotent}
$t$. The answer can be deduced from the Taylor Series at $0$ for the
function $K$, and the coefficients can be read off from Proposition
\ref{limit3}; alternatively, by the uniqueness of formal solutions
(Theorem \ref{formal}), they can be read off from the formal solution
we know already from Proposition \ref{exp-sol}. In any case, we get
for nilpotent $t$
\begin{equation}
K(t)= \delta (0) + t\cdot \Delta (\delta (0)) + \frac{t^2}{2!} \Delta
^2 (\delta (0)) + \ldots 
\label{formalheat}\end{equation}
the series being a finite sum, since $t$ is nilpotent. In particular,
for $d$ with $d^2 =0$, we have $K(d)= \delta (0) + d\cdot \Delta
(\delta (0))$, or since $\Delta = (\; )''$,
\begin{equation}K(d)= \delta (0) + d\cdot \delta (0)'' 
\label{littleheat}\end{equation}
In some sense, the motivation for our study of the heat equation 
in particular was to see how $\delta (0)$ evolves in nilpotent lapse $t$ of 
time and specially for $t=d$ with $d^2 =0$; the answer is
(\ref{littleheat}) (or more generally (\ref{formalheat})).

Being an extensive quantity, a distribution like (\ref{littleheat})
should be drawable.  In fact, it can be exhibited as a finite linear
combination of Dirac distributions $\delta (a)$ (= ``evaluate at $a$'').  This 
hinges on:

\begin{proposition}  Let $h^4 =0$.  Then
$$h^2 \cdot \delta (0)'' = \delta (-h) -2\delta (0) + \delta
(h).\label{sporer}$$
\end{proposition}
{\bf Proof.} It suffices to prove, for  an arbitrary test
function $\phi $, that
$h^2 \cdot \phi '' (0) = \phi (-h ) - 2\phi (0) + \phi (h)$; now just
Taylor expand the two outer terms in the sum on the right; the terms
of odd degree cancel, the terms of even degree (0 and two) give the
result.  (There is a similar result for higher derivatives of $\delta
(0)$: for $h^{n+1}=0$,
$$h^n \cdot \delta (0)^{(n)} = \sum _{i=0}^n (-1)^{i} (n,i) \delta (i\cdot h),$$
where $(n,i)$ denotes the binomial coefficient $n!/i!  (n-i)!$.
This hinges on some combinatorics with binomial coefficients, cf.\  
\cite{fe} p.  63, Problem 16).

\medskip
To make a ``drawing'' of $K(d)$ where $d^2 =0$, we assume that $d =
h^3$ for some $h$ with $h^4 =0$ (we shall not deal here with the question
whether this can always be done).  Then

$$K(d)= \delta (0) + d\cdot \delta (0)'' =
\delta (0) + h\cdot h^2 \cdot \delta (0)'' =
\delta (0) + h\cdot ((\delta (-h) -2\delta (0) + \delta (h))$$

\noindent using (\ref{littleheat}) and (\ref{sporer}).  The drawing one can
make of $\delta (x)$ (as for any discrete distribution), is a column
diagram: erect a column of heigth 1 at $x$.  The distribution above
then comes about by removing $2h$ units from the unit column at $0$,
and placing the small columns of heigth $h$ at $-h$ and $h$.  This is
the beginning of the {\em diffusion} of the Dirac distribution.
Several other ways of exhibiting $K(d)$ as linear combination of Dirac
distributions are also possible.

\medskip

Since ${\cal D}'(R^n )$ is a microlinear and Euclidean $R$-module,
and $\Delta : {\cal D}'(R^n )\to {\cal D}'(R^n )$ is linear,
 we may 
apply the general results of Propositions \ref{exp-sol} and 
\ref{exp-sol2} to conclude
that the formal solution of the heat equation $\dot{F}(t)=\Delta (F(t))$
 with initial value (the 
distribution) $\mu$, is
 the series
$$  \mu + t\Delta (\mu) + t^2/2! {\Delta}^2(\mu) + t^3/3!
{\Delta}^3(\mu) + \dots .$$
Similarly,
the formal solution of the wave equation $\ddot{F}(t)=\Delta (F(t))$ with 
initial value (the distribution) $\mu$, and initial speed the 
distribution $\nu$ is the series $$ \mu + t\nu  + t^2/2!  
{\Delta}(\mu) + t^3/3!
{\Delta}(\nu) + t^4 /4!  \Delta ^2 (\mu ) + t^5 /5!  \Delta ^2 (\nu )
\dots .$$

 Applying (in the one-variable case, say) these formulas to a test function 
 $\phi$ in the variable $x$ and to the distributions $\delta (0)$ and 
 $\delta'(0)$ we obtain the following Maclaurin series for the heat 
 equation 
 $$<F(t), \phi> = \phi (0) + t {\phi}^{('')}(0) + t^2/2!
{\phi}^{('''')}(0) +\dots .$$
Here $\dot{()}$ refers to the time derivative, whereas $()^{'}$ to the
space derivative $\partial/\partial x.$ The variable $x$ has been left
unexpressed.  There is a similar series for the wave equation:
$$<F(t), \phi> = \phi (0) + t \dot{\phi}(0) + t^2/2! 
{\phi}^{('')}(0) + t^3/3! {\dot{\phi}}^{('')}(0) + t^4/4! \phi^{''''}(0) +\dots .$$
\subsection{Simple Transport}

For the sake of completeness, we also consider the function $\delta : R\to {\bf D}' _c 
(R)$ given by $t\mapsto \delta (t)$, the Dirac distribution at $t\in
R$.  This is the ``fundamental solution'' for the equation for
``simple transport'', cf.\ e.g.\ \cite{str}.

\begin{proposition}The function $\delta$ is the solution for the
differential equation for ``simple transport",
$$\frac{d}{dt}(\delta ) = (\delta )'$$
with initial value $\delta (0)$.
\label{simple}\end{proposition}

{\bf Proof.} For any test function $\phi$,
$$\frac{d}{dt}<\delta (t), \phi > = \frac{d}{dt}\phi (t) = \phi ' (t) =
<\delta (t)' , \phi >.$$

\end{document}